\documentclass[11pt,a4paper]{article}

\usepackage{inputenc}
\usepackage{amsmath}
\usepackage{bm}
\usepackage{bbold}
\usepackage{amsthm}

\usepackage{hyperref}

\title{The Max-Plus Algebra Approach\\ in Modelling of Queueing Networks\thanks{Proceedings of the Summer Computer Simulation Conference, July 21-25, 1996, Portland, OR / Ed. by V. W. Ingalls, J. Cynamon, A. Saylor. The Society for Computer Simulation, San Diego, CA, 1996. P.~485-490.}}

\author{Nikolai K. Krivulin\thanks{Faculty of Mathematics and Mechanics, St.~Petersburg State University, 28 Universitetsky Ave., St.~Petersburg, 198504, Russia, 
nkk@math.spbu.ru}
}

\date{}

\newtheorem{theorem}{Theorem}
\newtheorem{lemma}[theorem]{Lemma}

\newtheorem{example}{Example}

\def\sumo_#1^#2{\setbox0=\hbox{$\displaystyle{\sum}$}
                \setbox1=\hbox{$\scriptstyle{#1}$}
                \setbox2=\hbox{$\scriptstyle{#2}$}
		\setbox3=\hbox{${}_{{}_\oplus}\mathsurround=0pt$}
		\dimen1=.5\wd1 \advance\dimen1 by-.5\wd0
		\ifdim\dimen1>0pt
		   \ifdim\dimen1>\wd3 \kern\wd3 \else\kern\dimen1\fi\fi
		\dimen2=.5\wd2 \advance\dimen2 by-.5\wd0
		\ifdim\dimen2>0pt
		   \ifdim\dimen2>\wd3 \kern\wd3 \else\kern\dimen2\fi\fi
		\mathop{{\sum}{}_{{}_\oplus}}_{\kern-\wd3 #1}^{\kern-\wd3 #2}}
\def\sumol_#1{\setbox0=\hbox{$\displaystyle{\sum}$}
             \setbox1=\hbox{$\scriptstyle{#1}$}
	     \setbox3=\hbox{${}_{{}_\oplus}\mathsurround=0pt$}
	     \dimen1=.5\wd1 \advance\dimen1 by-.5\wd0
	     \ifdim\dimen1>0pt
	        \ifdim\dimen1>\wd3 \kern\wd3 \else\kern\dimen1\fi\fi
	     \mathop{{\sum}{}_{{}_\oplus}}_{\kern-\wd3 #1}}

\begin{document}

\maketitle

\begin{abstract}
A class of queueing networks which consist of single-server fork-join nodes
with infinite buffers is examined to derive a representation of the network
dynamics in terms of max-plus algebra. For the networks, we present a common
dynamic state equation which relates the departure epochs of customers
from the network nodes in an explicit vector form determined by a state
transition matrix. We show how the matrix may be calculated from the service
time of customers in the general case, and give examples of matrices inherent
in particular networks.
\\

\end{abstract}

\section{Introduction}
We consider a class of queueing networks with single-server nodes and
customers of a single class. The server at each node is supplied with an
infinite buffer intended for customers waiting for service. There is, in
general, no restriction on the network topology; in particular, both open and
closed networks may be included in the class.

In addition to the ordinary service procedure, specific fork-join operations
\cite{Bacc89,Gree91} may be performed in each node of the network. In fact,
these operations allow customers (jobs, tasks) to be split into parts, and to
be merged into one, when circulating through the network.

Furthermore, we assume that there may be nodes, each distributing customers
among a group of its downstream nodes by a regular round routing mechanism.
Every node operating on the round routing basis passes on the consecutive
customers to distinct nodes being assigned to each new customer cyclically in
a fixed order.

Both the fork-join formalism and the above regular round routing scheme seem
to be useful in the description of dynamical processes in a variety of actual
systems, including production processes in manufacturing, transmission of
messages in communication networks, and parallel data processing in
multi-processor computer systems (see, e.g., \cite{Bacc89}).

In this paper, the networks are examined so as to represent their dynamics in
terms of max-plus algebra \cite{Cuni79,Olsd92}. The max-plus algebra approach
actually offers a quite compact and unified way of describing system dynamics,
which may provide a useful basis for both analytical study and simulation of
queueing systems. In particular, it has been shown in \cite{Kriv94,Kriv95} that
the evolution of both open and closed tandem queueing systems may be described
by the linear algebraic equation
\begin{equation}\label{ve-sde}
\bm{d}(k) = T(k) \otimes \bm{d}(k-1),
\end{equation}
where $  \bm{d}(k)  $ is a vector of departure epochs from the
queues, $  T(k)  $ is a matrix calculated from service times of customers,
and $  \otimes  $ is an operator which determines the matrix-vector
multiplication in the max-plus algebra.

The purpose of this paper is to demonstrate that the network dynamics also
allows of representation in the form of dynamic state equation (\ref{ve-sde}).
We start with preliminary max-plus algebra definitions and related results.
Furthermore, we give a general description of a queueing network model, which
is then refined to describe fork-join networks. It is shown how the dynamics
of the networks may be represented in the form of (\ref{ve-sde}).

The obtained representation is extended to describe the dynamics of tandem
queueing systems and a system with regular round routing. In fact, tandem
systems may be treated as trivial fork-join networks, and thus described in
the same way. To represent the system with round routing, which operates
differently than the other networks under examination, we first introduce an
equivalent fork-join network, and then get equation (\ref{ve-sde}).

\section{Algebraic Definitions and Results}
We start with a brief overview of basic facts about max-plus algebra, which we
will exploit in the representation of queueing network models in the
subsequent sections. Further details concerning the max-plus algebra theory
as well as its applications can be found in \cite{Cuni79,Cohe89,Olsd92}.

Max-plus algebra is normally defined (see, for example, \cite{Olsd92}) as the
system $  \langle \underline{\mathbb{R}}, \oplus, \otimes \rangle $, where
$  \underline{\mathbb{R}} = \mathbb{R} \cup \{\varepsilon\}  $ with
$  \varepsilon = -\infty $, and for any
$  x,y \in \underline{\mathbb{R}} $,
$$
x \oplus y = \max(x,y), \quad x \otimes y = x + y.
$$

Since the new operations $  \oplus  $ and $  \otimes  $ retain most of the
properties of the ordinary addition and multiplication, including
associativity, commutativity, and distributivity of multiplication over
addition, one can perform usual algebraic manipulations in the max-plus
algebra under the standard conventions regarding brackets and precedence of
$  \otimes  $ over $  \oplus $. However, the operation $  \oplus  $ is
idempotent; that is, for any $  x \in \underline{\mathbb{R}} $, one has
$  x \oplus x = x $.

There are the null and identity elements in this algebra, namely
$  \varepsilon  $ and $  0 $, which satisfy the evident conditions
$  x \oplus \varepsilon = \varepsilon \oplus x = x $, and
$  x \otimes 0 = 0 \otimes x = x $, for any
$  x \in \underline{\mathbb{R}} $. The absorption rule involving
$  x \otimes \varepsilon = \varepsilon \otimes x = \varepsilon  $ is also
true in the algebra.

The max-plus algebra of matrices is readily introduced in the regular way.
Specifically, for any $(n \times n)$-matrices $  X = (x_{ij})  $ and
$  Y = (y_{ij}) $, the entries of $  U = X \oplus Y  $ and
$  V = X \otimes Y  $ are calculated as
$$
u_{ij} = x_{ij} \oplus y_{ij}, \quad \mbox{and} \quad
v_{ij} = \sumo_{k=1}^{n} x_{ik} \otimes y_{kj},
$$
where $  \sum_{\oplus}  $ stands for the iterated operation $  \oplus $. As
the null element, the matrix $  \mathcal{E}  $ with all entries equal to
$  \varepsilon  $ is taken in the algebra; the matrix
$  E = \mathop\mathrm{diag}(0,\ldots,0)  $ with the off-diagonal entries set to
$  \varepsilon  $ presents the identity.

As is customary, for any square matrix $  X $, one can define $  X^{0} = E  $ and
$  X^{q} = X \otimes X^{q-1} = X^{q-1} \otimes X  $ for all
$  q=1,2,\ldots  $. Note, however, that idempotency in this algebra leads, in
particular, to the matrix identity
$$
(E \oplus X)^{q} = E \oplus X \oplus \cdots \oplus X^{q}.
$$

Many phenomena inherent in the max-plus algebra appear to be well explained in
terms of their graph interpretations \cite{Cuni79,Cohe89,Olsd92}. To
illustrate, consider an $(n \times n)$-matrix $  X  $ with entries
$  x_{ij} \in \underline{\mathbb{R}} $, and note that it can be treated as
the adjacency matrix of an oriented graph with $  n  $ nodes, provided each
entry $  x_{ij} \neq \varepsilon  $ implies the existence of the arc
$  (i,j)  $ in the graph, whereas $  x_{ij} = \varepsilon  $ does the lack
of the arc.

Let us calculate the matrix $  X^{2} = X \otimes X $, and denote its entries
by $  x_{ij}^{(2)} $. Clearly, we have $  x_{ij}^{(2)} \neq \varepsilon  $
if and only if there exists at least one path from node $  i  $ to node
$  j  $ in the graph, which consists of two arcs. Moreover, for any integer
$  q > 0 $, the matrix $  X^{q}  $ has the entry
$  x_{ij}^{(q)} \neq \varepsilon  $ only if there is a path with the length
$  q  $ from $  i  $ to $  j $.

Suppose that the graph associated with $  X  $ is acyclic. It is clear that
we will have $  X^{q} = \mathcal{E}  $ for all $  q > p $, where $  p  $ is
the length of the longest path in the graph. Assume now the graph not to be
acyclic, and then consider any one of its circuits. Since it is possible to
construct a cyclic path of any length, which lies along the circuit, we
conclude that $  X^{q} \neq \mathcal{E}  $ for all $  q=1,2,\ldots $

Finally, we consider the implicit equation in the vector
$  \bm{x} = (x_{1},\ldots,x_{n})^{T} $,
\begin{equation}\label{e-I}
\bm{x} = U \otimes \bm{x}
                     \oplus \bm{v},
\end{equation}
where $  U = (u_{ij})  $ and
$  \bm{v} = (v_{1},\ldots,v_{n})^{T}  $ are respectively given
$(n\times n)$-matrix and $n$-vector. This equation actually plays a large
role in max-plus algebra representations of dynamical systems including
systems of queues \cite{Cohe89,Olsd92,Kriv95}. The next lemma offers
particular conditions for (\ref{e-I}) to be solvable, and it shows how the
solution may be calculated. One can find a detailed analysis of (\ref{e-I}) in
the general case in \cite{Cohe89}.

\begin{lemma}\label{l-1}
Suppose that the entries of the matrix $  U  $ and the vector
$  \bm{v}  $ are either positive or equal to
$  \varepsilon $. Then equation {\rm (\ref{e-I})} has the unique bounded
solution $  \bm{x}  $ if and only if the graph associated with
the matrix $  U  $ is acyclic. Provided that the solution exists, it is
given by
$$
\bm{x} = (E \oplus U)^{p} \otimes \bm{v},
$$
where $  p  $ is the length of the longest path in the graph.
\end{lemma}

The proof of the lemma may be furnished based on the above graph
interpretation and idempotency of $  \oplus $.

\section{A General Network Description}
In this section, we introduce some general notations and give related
definitions, which are common for all particular networks examined below.
In fact, we consider a network with $  n  $ single-server nodes and
customers of a single class. The network topology is described by an oriented
graph $  \mathcal{G} = ({\bf N}, {\bf A}) $, where
$  {\bf N} = \{1,\ldots,n\}  $ represents the nodes, and
$  {\bf A} = \{(i,j)\} \subset {\bf N} \times {\bf N}  $ does the arcs
determining the transition routes of customers.

For every node $  i \in {\bf N} $, we introduce the set of its predecessors
$  {\bf P}(i) = \{ j | \, (j,i) \in {\bf A} \}  $ and the set of its
successors $  {\bf S}(i) = \{ j | \, (i,j) \in {\bf A} \} $. In specific
cases, there may be one of the conditions $  {\bf P}(i) = \emptyset  $ and
$  {\bf S}(i) = \emptyset  $ encountered. Each node
$  i  $ with $  {\bf P}(i) = \emptyset  $ is assumed to represent an
infinite external arrival stream of customers; provided that
$  {\bf S}(i) = \emptyset $, it is considered as an output node intended to
release customers from the network.

Each node $  i \in {\bf N}  $ includes a server and its buffer with infinite
capacity, which together present a single-server queue operating under the
first-come, first-served (FCFS) queueing discipline. At the initial time, the
server at each node $  i  $ is assumed to be free of customers, whereas in
its buffer, there may be $  r_{i} $, $  0 \leq r_{i} \leq \infty $,
customers waiting for service. The value $  r_{i} = \infty  $ is set for
every node $  i  $ with $  {\bf P}(i) = \emptyset $, which represents an
external arrival stream of customers.

For the queue at node $  i $, we denote the $k$th arrival and departure
epochs respectively as $  a_{i}(k)  $ and $  d_{i}(k) $. Furthermore, the
service time of the $k$th customer at server $  i  $ is indicated by
$  \tau_{ik} $. We assume that $  \tau_{ik} > 0  $ are given parameters for
all $  i=1,\ldots,n $, and $  k=1,2,\ldots $, while $  a_{i}(k)  $ and
$  d_{i}(k)  $ are considered as unknown state variables. With the condition
that the network starts operating at time zero, it is convenient to set
$  d_{i}(0) \equiv 0 $, and $  d_{i}(k) \equiv \varepsilon  $ for all
$  k < 0 $, $  i=1,\ldots,n $.

It is easy to set up an equation which relates the system state variables. In
fact, the dynamics of any single-server node $  i  $ with an infinite
buffer, operating on the FCFS basis, is described as \cite{Kriv95}
\begin{equation}\label{se-d}
d_{i}(k) = \tau_{ik} \otimes a_{i}(k) \oplus \tau_{ik} \otimes d_{i}(k-1).
\end{equation}
With the system state vectors
$$
\bm{a}(k)
= \left(
    \begin{array}{c}
      a_{1}(k) \\
      \vdots   \\
      a_{n}(k)
    \end{array}
  \right),
\quad
\bm{d}(k)
= \left(
    \begin{array}{c}
      d_{1}(k) \\
      \vdots   \\
      d_{n}(k)
    \end{array}
  \right),
$$
and the diagonal matrix
$$
\mathcal{T}_{k}
= \left(
    \begin{array}{ccc}
      \tau_{1k}   &        & \varepsilon \\
                  & \ddots &             \\
      \varepsilon &        & \tau_{nk}
    \end{array}
  \right),
$$
we may rewrite equation (\ref{se-d}) in a vector form, as
\begin{equation}\label{ve-d}
\bm{d}(k) = \mathcal{T}_{k} \otimes \bm{a}(k)
                       \oplus \mathcal{T}_{k} \otimes \bm{d}(k-1).
\end{equation}

Clearly, to represent the dynamics of a network completely, equation
(\ref{ve-d}) should be supplemented with that determining the vector of
arrival epochs, $  \bm{a}(k) $. The latter equation may differ
for distinct networks according to their operation features and topology. We
will give appropriate equations for $  \bm{a}(k)  $ inherent
in particular networks, as well as related representations of the entire
network dynamics in the subsequent sections.

\section{Fork-Join Queueing Networks}
The purpose of this section is to derive an algebraic representation of the
dynamics of fork-join networks which present a quite general class of queueing
network models. Since we do not impose any limitation on the network topology,
the models under study may be considered as an extension of acyclic fork-join
queueing networks investigated in \cite{Bacc89}.

The distinctive feature of any fork-join network is that, in addition to the
usual service procedure, special join and fork operations are performed in its
nodes, respectively before and after service. The join operation is actually
thought to cause each customer which comes into node $  i $, not to enter the
buffer at the server but to wait until at least one customer from every node
$  j \in {\bf P}(i)  $ arrives. As soon as these customers arrive, they,
taken one from each preceding node, are united to be treated as being one
customer which then enters the buffer to become a new member of the queue.

The fork operation at node $  i  $ is initiated every time the service of a
customer is completed; it consists in giving rise to several new customers
instead of the original one. As many new customers appear in node $  i  $ as
there are succeeding nodes included in the set $  {\bf S}(i) $. These
customers simultaneously depart the node, each being passed to separate node
$  j \in {\bf S}(i) $. We assume that the execution of fork-join operations
when appropriate customers are available, as well as the transition of
customers within and between nodes require no time.

As it immediately follows from the above description of the fork-join
operations, the $k$th arrival epoch into the queue at node $  i  $ is
represented as (also, see \cite{Bacc89,Gree91})
\begin{equation}\label{se-a}
a_{i}(k) = \left\{
	     \begin{array}{ll}
	       \displaystyle{\sumol_{j\in \mbox{\scriptsize\bf P}(i)}}
                                                               d_{j}(k-r_{i}),
                               & \mbox{if $  {\bf P}(i) \neq \emptyset $}, \\
	       \quad\varepsilon, & \mbox{if $  {\bf P}(i) = \emptyset $}.
	     \end{array}
	   \right.
\end{equation}

In order to get this equation in a vector form, we first define
$  M = \max \{ r_{i} | \, r_{i} < \infty, \, i=1,\ldots,n \} $. Now we
may rewrite (\ref{se-a}) as
$$
a_{i}(k) = \sumo_{m=0}^{M} \sumo_{j=1}^{n} g_{ji}^{m} \otimes d_{j}(k-m),
$$
where the numbers $  g_{ij}^{m}  $ are determined by the condition
\begin{equation}\label{e-g}
g_{ij}^{m} = \left\{
	       \begin{array}{ll}
                0, & \mbox{if $  i \in {\bf P}(j)  $ and $  m=r_{j} $}, \\
                \varepsilon, & \mbox{otherwise}.
               \end{array}
	     \right.
\end{equation}

Let us introduce the matrices $  G_{m} = \left(g_{ij}^{m}\right)  $ for each
$  m=0,1,\ldots,M $, and note that $  G_{m}  $ actually presents an
adjacency matrix of the partial graph
$  \mathcal{G}_{m} = ({\bf N},{\bf A}_{m}) $, where
$  {\bf A}_{m} = \{(i,j)| \, i \in {\bf P}(j), \, r_{j} = m \} $. With these
matrices, equation (\ref{se-a}) may be written in the vector form
\begin{equation}\label{ve-a}
\bm{a}(k) = \sumo_{m=0}^{M} G_{m}^{T}
                                            \otimes \bm{d}(k-m),
\end{equation}
where $  G_{m}^{T}  $ denotes the transpose of the matrix $  G_{m} $.

Furthermore, by combining equations (\ref{ve-d}) and (\ref{ve-a}), we arrive
at the equation
\begin{eqnarray}
\bm{d}(k)  & =  & \mathcal{T}_{k} \otimes G_{0}^{T}
                                            \otimes \bm{d}(k)
            \oplus \mathcal{T}_{k} \otimes \bm{d}(k-1) \nonumber \\
& & \mbox{} \oplus \mathcal{T}_{k} \otimes \sumo_{m=1}^{M} G_{m}^{T}
                              \otimes \bm{d}(k-m) \label{ve-ex}.
\end{eqnarray}
Clearly, it is actually an implicit equation in $  \bm{d}(k) $,
which has the form of (\ref{e-I}), with 
$  U = \mathcal{T}_{k} \otimes G_{0}^{T} $. Taking into account that the matrix 
$  \mathcal{T}_{k}  $ is diagonal, one can apply Lemma~\ref{l-1} to prove the 
following statement.
\begin{theorem}\label{t-2}
Suppose that in the fork-join network model, the graph $  \mathcal{G}_{0}  $
associated with the matrix $  G_{0}  $ is acyclic. Then equation
{\rm (\ref{ve-ex})} can be solved to produce the explicit dynamic state
equation
\begin{equation}\label{ve-sde1}
\bm{d}(k) = \sumo_{m=1}^{M} T_{m}(k)
                                            \otimes \bm{d}(k-m),
\end{equation}
with the state transition matrices
\begin{eqnarray*}
T_{1}(k) & = & (E \oplus \mathcal{T}_{k} \otimes G_{0}^{T})^{p}
                                   \otimes \mathcal{T}_{k}
                                              \otimes (E \oplus G_{1}^{T}), \\
T_{m}(k) & = & (E \oplus \mathcal{T}_{k} \otimes G_{0}^{T})^{p}
                                    \otimes \mathcal{T}_{k} \otimes G_{m}^{T}, \\
         &   & \mbox{} m=2,\ldots,M,
\end{eqnarray*}
where $  p  $ is the length of the longest path in $  \mathcal{G}_{0} $.
\end{theorem}

Finally, with the extended state vector
$$
\widehat{\bm{d}}(k) = \left(
                                      \begin{array}{l}
                                        \bm{d}(k) \\
                                        \bm{d}(k-1) \\
                                        \vdots \\
                                        \bm{d}(k-M+1)
                                      \end{array}
                                    \right),
$$
one can bring equation (\ref{ve-sde1}) into the form of (\ref{ve-sde}):
$$
\widehat{\bm{d}}(k)
                 = \widehat{T}(k) \otimes \widehat{\bm{d}}(k-1),
$$
where the state transition matrix is defined as
$$
\widehat{T}(k) = \left(
		   \begin{array}{ccccc}
                     T_{1}(k) & T_{2}(k) & \cdots & \cdots &T_{M}(k)  \\
		     E        & \mathcal{E} & \cdots & \cdots & \mathcal{E} \\
		              & \ddots   & \ddots &        & \vdots   \\
		              &          & \ddots & \ddots & \vdots   \\
		     \mathcal{E} &          &        & E      & \mathcal{E}
		   \end{array}
		 \right).
$$

To conclude this section, we present an example which shows how the dynamics
of a particular fork-join network is described based on the above
representation.
\begin{example} We consider a network with $  n=5  $ nodes, depicted in
Fig.~\ref{f-FJN}. The initial numbers of customers in the network nodes are
determined as follows: $  r_{1}=\infty $, $  r_{2}=r_{4}=0 $, and
$  r_{3}=r_{5}=1 $.
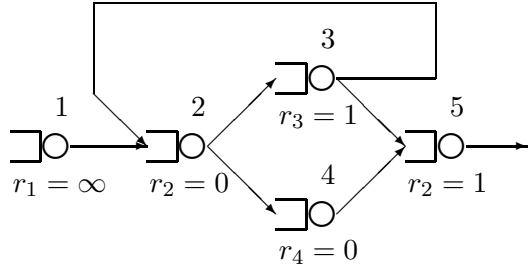
\begin{figure}[hhh]
\begin{center}
\begin{picture}(70,35)
\thicklines

\put(5.5,19){$1$}
\put(0,9){$r_{1}=\infty$}
\put(0,17){\line(1,0){4}}
\put(0,13){\line(1,0){4}}
\put(4,17){\line(0,-1){4}}
\put(6,15){\circle{3}}

\put(8,15){\thinlines\vector(1,0){10}}

\put(23.5,19){$2$}
\put(18,9){$r_{2}=0$}
\put(18,17){\line(1,0){4}}
\put(18,13){\line(1,0){4}}
\put(22,17){\line(0,-1){4}}
\put(24,15){\circle{3}}

\put(26,15){\thinlines\vector(1,1){9}}
\put(26,15){\thinlines\vector(1,-1){9}}

\put(40.5,28){$3$}
\put(35,18){$r_{3}=1$}
\put(35,26){\line(1,0){4}}
\put(35,22){\line(1,0){4}}
\put(39,26){\line(0,-1){4}}
\put(41,24){\circle{3}}

\put(43,24){\thinlines\line(1,0){13}}
\put(56,24){\thinlines\line(0,1){10}}
\put(56,34){\thinlines\line(-1,0){45}}
\put(11,34){\thinlines\line(0,-1){12}}
\put(11,22){\thinlines\vector(1,-1){7}}
\put(43,24){\thinlines\vector(1,-1){9}}

\put(40.5,10){$4$}
\put(35,0){$r_{4}=0$}
\put(35,8){\line(1,0){4}}
\put(35,4){\line(1,0){4}}
\put(39,8){\line(0,-1){4}}
\put(41,6){\circle{3}}
\put(43,6){\thinlines\vector(1,1){9}}

\put(57.5,19){$5$}
\put(52,9){$r_{2}=1$}
\put(52,17){\line(1,0){4}}
\put(52,13){\line(1,0){4}}
\put(56,17){\line(0,-1){4}}
\put(58,15){\circle{3}}
\put(60,15){\thinlines\vector(1,0){8}}

\end{picture}
\end{center}
\caption{A fork-join queueing network.}\label{f-FJN}
\end{figure}

First note that for the network, we have $  M=1 $. In this case, provided
that explicit representation (\ref{ve-sde1}) exists, it is just written in the
form of (\ref{ve-sde}):
$$
\bm{d}(k) = T(k) \otimes \bm{d}(k-1),
$$
with
$  T(k) = T_{1}(k) = (E \oplus \mathcal{T}_{k} \otimes G_{0}^{T})^{p}
 \otimes \mathcal{T}_{k} \otimes (E \oplus G_{1}^{T}) $.

Furthermore, by applying (\ref{e-g}), we may calculate
{\setlength{\arraycolsep}{1.5mm}
$$
G_{0}
= \left(
    \begin{array}{ccccc}
      \varepsilon & 0           & \varepsilon & \varepsilon & \varepsilon \\
      \varepsilon & \varepsilon & \varepsilon & 0           & \varepsilon \\
      \varepsilon & 0           & \varepsilon & \varepsilon & \varepsilon \\
      \varepsilon & \varepsilon & \varepsilon & \varepsilon & \varepsilon \\
      \varepsilon & \varepsilon & \varepsilon & \varepsilon & \varepsilon
    \end{array}
  \right), \quad
G_{1}
= \left(
    \begin{array}{ccccc}
      \varepsilon & \varepsilon & \varepsilon & \varepsilon & \varepsilon \\
      \varepsilon & \varepsilon & 0           & \varepsilon & \varepsilon \\
      \varepsilon & \varepsilon & \varepsilon & \varepsilon & 0           \\
      \varepsilon & \varepsilon & \varepsilon & \varepsilon & 0           \\
      \varepsilon & \varepsilon & \varepsilon & \varepsilon & \varepsilon
    \end{array}
  \right).
$$
}%

Since the graph associated with the matrix $  G_{0}  $ is acyclic, with the
length of its longest path $  p=2 $, we may really describe the network
dynamics through equation (\ref{ve-sde}). Finally, simple algebraic
manipulations give
{\setlength{\arraycolsep}{1.0mm}
\begin{eqnarray*}
\lefteqn{T(k) = (E \oplus \mathcal{T}_{k} \otimes G_{0}^{T})^{2}
                       \otimes \mathcal{T}_{k} \otimes (E \oplus G_{1}^{T}) =} \\
&& \left(
    \begin{array}{ccccc}
      \tau_{1k}   & \varepsilon & \varepsilon & \varepsilon & \varepsilon \\
      \tau_{1k} \otimes \tau_{2k} &
      \tau_{2k} \otimes \tau_{3k} &
                  \tau_{2k} \otimes \tau_{3k} & \varepsilon & \varepsilon \\
      \varepsilon & \tau_{3k}   & \tau_{3k}   & \varepsilon & \varepsilon \\
      \tau_{1k} \otimes \tau_{2k} \otimes \tau_{4k} &
      \tau_{2k} \otimes \tau_{3k} \otimes \tau_{4k} &
      \tau_{2k} \otimes \tau_{3k} \otimes \tau_{4k} &
                                                  \tau_{4k} & \varepsilon \\
      \varepsilon & \varepsilon & \tau_{5k}   & \tau_{5k}   & \tau_{5k}
    \end{array}
  \right).
\end{eqnarray*}
}%
\end{example}

\section{Tandem Queues}
Tandem queueing systems present networks with the simplest topology determined
by graphs which include only the nodes with no more than one incoming and
outgoing arcs. Although no fork and join operations are actually performed in
the systems, yet they may be treated as trivial fork-join networks and thus
described using the representation proposed in the previous section.

\begin{figure}[hhh]
\begin{center}
\begin{picture}(65,15)
\thicklines

\put(5.5,10){$1$}
\put(0,0){$r_{1}=\infty$}
\put(0,8){\line(1,0){4}}
\put(0,4){\line(1,0){4}}
\put(4,8){\line(0,-1){4}}
\put(6,6){\circle{3}}
\put(8,6){\thinlines\vector(1,0){8}}

\put(21.5,10){$2$}
\put(16,0){$r_{2}=0$}
\put(16,8){\line(1,0){4}}
\put(16,4){\line(1,0){4}}
\put(20,8){\line(0,-1){4}}
\put(22,6){\circle{3}}
\put(24,6){\thinlines\vector(1,0){8}}

\multiput(34,6)(2,0){3}{\circle*{1}}
\put(40,6){\thinlines\vector(1,0){8}}

\put(53.5,10){$n$}
\put(48,0){$r_{n}=0$}
\put(48,8){\line(1,0){4}}
\put(48,4){\line(1,0){4}}
\put(52,8){\line(0,-1){4}}
\put(54,6){\circle{3}}
\put(56,6){\thinlines\vector(1,0){8}}

\end{picture}
\end{center}
\caption{Open tandem queues.}\label{f-OTQ}

\end{figure}
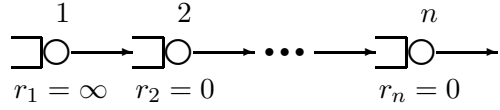

Let us first consider a series of $  n  $ single-server queues, depicted in
Fig.~\ref{f-OTQ}. In this open tandem system, the queue labelled with
$  1  $ is assigned to represent an infinite external arrival stream of
customers, with $  r_{1}=\infty $. The buffers of servers $  2  $ to
$  n  $ are assumed to be empty at the initial time, so we set
$  r_{i}=0  $ for all $  i=2,\ldots,n $.

As one can see, equation (\ref{se-a}) is now reduced to
$$
a_{i}(k) = \left\{
             \begin{array}{ll}
               \varepsilon,   & \mbox{if $  i=1 $}, \\
               d_{i-1}(k),    & \mbox{if $  i\neq 1 $},
             \end{array}
           \right.
$$

Since $  M=0 $, there is only one matrix
$$
G_{0} = \left(
          \begin{array}{cccc}
            \varepsilon & 0      &        & \varepsilon \\
            \vdots      & \ddots & \ddots & \\
            \vdots      &        & \ddots & 0 \\
            \varepsilon & \cdots & \cdots & \varepsilon
          \end{array}
        \right)
$$
included in representation (\ref{ve-a}). It just presents an adjacency matrix
of the graph describing the topology of the open tandem system. It is clear
that this graph is acyclic; the length of its longest path $  p = n-1 $.

Finally, with $  G_{1} = \mathcal{E} $, equation (\ref{ve-sde1}) becomes
$$
\bm{d}(k)
= (E \oplus \mathcal{T}_{k} \otimes G_{0}^{T})^{n-1}
                       \otimes \mathcal{T}_{k} \otimes \bm{d}(k-1),
$$
which may be readily rewritten in the form of (\ref{ve-sde}) with
\begin{eqnarray*}
\lefteqn{T(k) = (E \oplus \mathcal{T}_{k} \otimes G_{0}^{T})^{n-1}
                                                      \otimes \mathcal{T}_{k}} \\
&&
= \left(
  \begin{array}{llcc}
    \tau_{1k}                   & \varepsilon & \cdots & \varepsilon \\
    \tau_{1k} \otimes \tau_{2k} & \tau_{2k}   &        & \varepsilon \\
    \vdots                      & \vdots      &        & \\
    \tau_{1k} \otimes \cdots \otimes \tau_{nk} &
    \tau_{2k} \otimes \cdots \otimes \tau_{nk} & \cdots & \tau_{nk}
  \end{array}
\right).
\end{eqnarray*}

\begin{figure}[hhh]
\begin{center}
\begin{picture}(75,20)
\thicklines

\put(13.5,14){$1$}
\put(9,4){$r_{1}$}
\put(0,10){\thinlines\vector(1,0){8}}
\put(8,12){\line(1,0){4}}
\put(8,8){\line(1,0){4}}
\put(12,12){\line(0,-1){4}}
\put(14,10){\circle{3}}
\put(16,10){\thinlines\vector(1,0){8}}

\put(29.5,14){$2$}
\put(25,4){$r_{2}$}
\put(24,12){\line(1,0){4}}
\put(24,8){\line(1,0){4}}
\put(28,12){\line(0,-1){4}}
\put(30,10){\circle{3}}
\put(32,10){\thinlines\vector(1,0){8}}

\multiput(42,10)(2,0){3}{\circle*{1}}
\put(48,10){\thinlines\vector(1,0){8}}

\put(61.5,14){$n$}
\put(57,4){$r_{n}$}
\put(56,12){\line(1,0){4}}
\put(56,8){\line(1,0){4}}
\put(60,12){\line(0,-1){4}}
\put(62,10){\circle{3}}
\put(64,10){\thinlines\vector(1,0){8}}

\put(0,10){\thinlines\line(0,-1){10}}
\put(72,10){\thinlines\line(0,-1){10}}
\put(0,0){\thinlines\line(1,0){72}}

\end{picture}
\end{center}
\caption{A closed tandem queueing system.}\label{f-CTQ}
\end{figure}
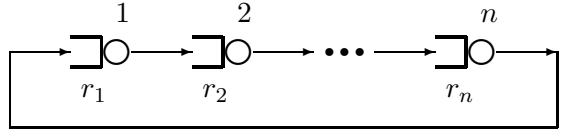

We now turn to a brief discussion of closed tandem systems. In the system
shown in Fig.~\ref{f-CTQ}, the customers pass through the queues consecutively
so as to get service at each server. After their service at queue $  n $, the
customers return to the $1$st queue for a new cycle of service.

There is a finite number of customers in the system at the initial time; that
is, we have $  0 \leq r_{i} < \infty  $ for all $  i=1,\ldots, n $. Note
that the numbers $  r_{i}  $ may all not equal $  0  $ since at least
one customer is to be present.

It is easy to understand that for the closed tandem system, equation
(\ref{se-a}) takes the form
$$
a_{i}(k) = \left\{
	     \begin{array}{ll}
                d_{n}(k-r_{1}),   & \mbox{if $  i=1 $}, \\
                d_{i-1}(k-r_{i}), & \mbox{if $  i\neq 1 $},
	     \end{array}
           \right.
$$
whereas equations (\ref{ve-a}) and (\ref{ve-ex}) remain unchanged. Since
there is at least one customer in the system, the graph $  \mathcal{G}_{0}  $
associated with the matrix $  G_{0}  $ cannot coincide with the system graph
$  \mathcal{G} $, and so is acyclic. As a consequence, the conditions of
Theorem~\ref{t-2} will be satisfied. We therefore conclude that the dynamics
of the system is described by equation (\ref{ve-sde1}) and thus by
(\ref{ve-sde}).

\begin{example} Let us suppose that in the system, $  r_{i}=1  $ for all
$  i=1,\ldots,n $. Then we have $  G_{0} = \mathcal{E}  $ and
$$
G_{1} = \left(
          \begin{array}{cccc}
            \varepsilon & 0           &             & \varepsilon \\
            \vdots      & \ddots      & \ddots      &             \\
            \varepsilon & \varepsilon & \ddots      & 0           \\
            0           & \varepsilon & \cdots      & \varepsilon
          \end{array}
        \right).
$$

Clearly, equation (\ref{ve-sde1}) is now written as
$$
\bm{d}(k)
= \mathcal{T}_{k} \otimes (E \oplus G_{1}^{T}) \otimes \bm{d}(k-1);
$$
that is, in the form of (\ref{ve-sde}) with the matrix
$$
T(k) = \mathcal{T}_{k} \otimes (E \oplus G_{1}^{T})
=
\left(
  \begin{array}{ccccc}
    \tau_{1k}   & \varepsilon & \cdots & \varepsilon & \tau_{1k}   \\
    \tau_{2k}   & \tau_{2k}   &        & \varepsilon & \varepsilon \\
                & \ddots      & \ddots &             &             \\
                &             & \ddots & \ddots      &             \\
    \varepsilon & \varepsilon &        & \tau_{nk}   & \tau_{nk}   \\
  \end{array}
\right).
$$
\end{example}

\section{A System with Round Routing}
\begin{figure}[hhh]
\begin{center}
\begin{picture}(45,55)
\thicklines

\put(5.5,30){$0$}
\put(-1,19){$r_{0}=\infty$}
\put(0,28){\line(1,0){4}}
\put(0,24){\line(1,0){4}}
\put(4,28){\line(0,-1){4}}
\put(6,26){\circle{3}}

\put(8,26){\thinlines\vector(1,1){20}}
\put(8,26){\thinlines\vector(4,1){20}}
\put(8,26){\thinlines\vector(1,-1){20}}

\put(33.5,50){$1$}
\put(28,40){$r_{1}=0$}
\put(28,48){\line(1,0){4}}
\put(28,44){\line(1,0){4}}
\put(32,48){\line(0,-1){4}}
\put(34,46){\circle{3}}
\put(36,46){\thinlines\vector(1,0){8}}

\put(33.5,35){$2$}
\put(28,25){$r_{2}=0$}
\put(28,33){\line(1,0){4}}
\put(28,29){\line(1,0){4}}
\put(32,33){\line(0,-1){4}}
\put(34,31){\circle{3}}
\put(36,31){\thinlines\vector(1,0){8}}

\multiput(34,20)(0,-2){3}{\circle*{1}}

\put(34,10){$l$}
\put(28,0){$r_{l}=0$}
\put(28,8){\line(1,0){4}}
\put(28,4){\line(1,0){4}}
\put(32,8){\line(0,-1){4}}
\put(34,6){\circle{3}}
\put(36,6){\thinlines\vector(1,0){8}}

\end{picture}
\end{center}
\caption{A queueing system with round routing.}\label{f-SRR}
\end{figure}
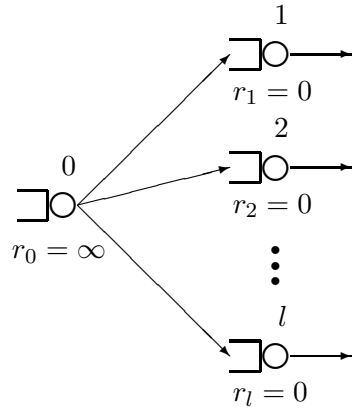

We consider an open system depicted in Fig.~\ref{f-SRR}, which consists of
$  n=l+1  $ queues labeled with $  0,1,\ldots,l $; queue $  0  $ is
thought to represent an external arrival stream of customers. Each incoming
customer has to go to one of the other queues, being chosen by a regular round
routing mechanism. Upon completion of his service at the latter queue, the
customer leaves the system.

The routing mechanism working in the system requires that the customer which
is the first to depart queue $  0 $, go to the $1$st queue, the second
customer do to the $2$nd queue, and so on. After the $l$th customer who is
transferred to queue $  l $, the next $(l+1)$st customer goes to the $1$st
queue once again, and the procedure is further repeated round and round.

\begin{figure}[hhh]
\begin{center}
\begin{picture}(80,65)
\thicklines

\put(10,38){$l+1$}
\put(4,28){$r_{l+1}=1$}
\put(0,34){\thinlines\vector(1,0){8}}
\put(8,36){\line(1,0){4}}
\put(8,32){\line(1,0){4}}
\put(12,36){\line(0,-1){4}}
\put(14,34){\circle{3}}
\put(16,34){\thinlines\line(1,1){22}}
\put(16,34){\thinlines\vector(1,-1){6}}

\put(67,60){$1$}
\put(61,50){$r_{1}=0$}
\put(38,56){\thinlines\vector(1,0){24}}
\put(62,58){\line(1,0){4}}
\put(62,54){\line(1,0){4}}
\put(66,58){\line(0,-1){4}}
\put(68,56){\circle{3}}
\put(70,56){\thinlines\vector(1,0){8}}

\put(24,32){$l+2$}
\put(18,22){$r_{l+2}=0$}
\put(22,30){\line(1,0){4}}
\put(22,26){\line(1,0){4}}
\put(26,30){\line(0,-1){4}}
\put(28,28){\circle{3}}
\put(30,28){\thinlines\line(1,1){12}}
\put(30,28){\thinlines\vector(1,-1){6}}

\put(67,44){$2$}
\put(60,34){$r_{2}=0$}
\put(42,40){\thinlines\vector(1,0){20}}
\put(62,42){\line(1,0){4}}
\put(62,38){\line(1,0){4}}
\put(66,42){\line(0,-1){4}}
\put(68,40){\circle{3}}
\put(70,40){\thinlines\vector(1,0){8}}

\multiput(37.5,20.5)(1.5,-1.5){3}{\circle*{1}}
\put(42,16){\thinlines\vector(1,-1){6}}
\multiput(68,30)(0,-2){3}{\circle*{1}}

\put(52,14){$2l$}
\put(46,4){$r_{2l}=0$}
\put(48,12){\line(1,0){4}}
\put(48,8){\line(1,0){4}}
\put(52,12){\line(0,-1){4}}
\put(54,10){\circle{3}}
\put(56,10){\thinlines\vector(1,1){6}}
\put(56,10){\thinlines\line(1,-1){6}}

\put(67,20){$l$}
\put(61,10){$r_{l}=0$}
\put(62,18){\line(1,0){4}}
\put(62,14){\line(1,0){4}}
\put(66,18){\line(0,-1){4}}
\put(68,16){\circle{3}}
\put(70,16){\thinlines\vector(1,0){8}}

\put(62,4){\thinlines\line(0,-1){4}}
\put(62,0){\thinlines\line(-1,0){62}}
\put(0,0){\thinlines\line(0,1){34}}
\put(0,34){\thinlines\vector(1,0){8}}

\end{picture}
\end{center}
\caption{An equivalent fork-join network.}\label{f-EFJ}
\end{figure}
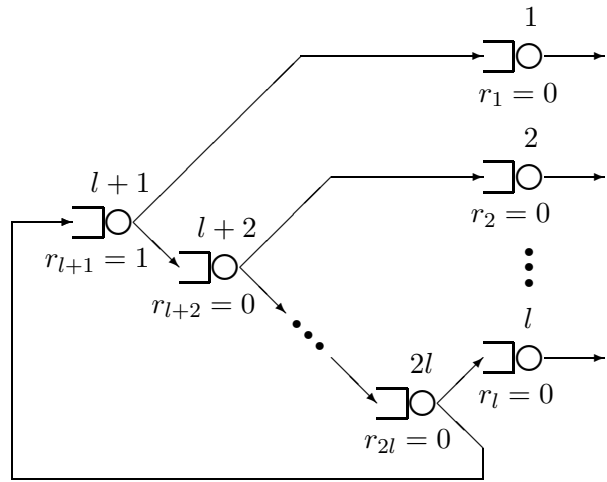

To represent the dynamics of the system, let us first replace it by an
equivalent fork-join network. An appropriate network, which is actually
obtained by substituting new nodes $  l+1  $ to $  2l  $ for queue
$  0 $, is shown in Fig.~\ref{f-EFJ}.

It is easy to see that, with the $k$th service time at each new node $  i $,
$  i=l+1,\ldots,2l $, defined as
$$
\tau_{ik} = \tau_{0\, lk-2l+i},
$$
the fork-join network behaves in much the same way as the original system with
regular round routing. Specifically, the overall operation of nodes
$  1  $ to $  l  $ in these queueing systems is identical, and thus both
systems produce the same output stream of customers.

The dynamics of the fork-join network is described through equations
(\ref{ve-d}) and (\ref{ve-a}) with $  n=2l  $ and $  M=1 $. Taking into
account that equation (\ref{se-a}) now becomes
$$
a_{i}(k) = \left\{
             \begin{array}{ll}
               d_{l+i}(k),   & \mbox{if $  i=1,\ldots, l $}, \\
               d_{2l}(k-1),  & \mbox{if $  i=l+1 $}, \\
               d_{i-1}(k),   & \mbox{if $  i=l+2,\ldots,n $}, \\
             \end{array}
           \right.
$$
we may represent the matrices $  G_{0}  $ and $  G_{1}  $ as
$$
G_{0} = \left(
          \begin{array}{cc}
            \mathcal{E}_{(l\times l)} & \mathcal{E}_{(l\times l)} \\
            E_{(l\times l)}        & F_{(l\times l)}
          \end{array}
        \right),
\quad
G_{1} = \left(
          \begin{array}{cc}
            \mathcal{E}_{(l\times l)} & \mathcal{E}_{(l\times l)} \\
            \mathcal{E}_{(l\times l)} & H_{(l\times l)}
          \end{array}
        \right),
$$
where $  \mathcal{E}_{(l\times l)}  $ and $  E_{(l\times l)}  $ respectively
denote the null and unit $(l\times l)$-matrices,
{\setlength{\arraycolsep}{1.2mm}
$$
F_{(l\times l)}
= \left(
    \begin{array}{cccc}
      \varepsilon & 0      &        & \varepsilon \\
      \vdots      & \ddots & \ddots &             \\
      \vdots      &        & \ddots & 0           \\
      \varepsilon & \cdots & \cdots & \varepsilon
    \end{array}
  \right),
\quad
H_{(l\times l)}
= \left(
    \begin{array}{cccc}
      \varepsilon & \cdots      & \cdots & \varepsilon \\
      \vdots      & \ddots      &        & \vdots      \\
      \varepsilon &             & \ddots & \vdots      \\
      0           & \varepsilon & \cdots & \varepsilon
    \end{array}
  \right).
$$
}%

As it is easy to verify, the graph associated with the matrix $  G_{0}  $ is
acyclic, with the length of its longest path $  p=l $. In this case, one can
apply Theorem~\ref{t-2} so as to obtain equation (\ref{ve-sde}) with the state
transition matrix
$$
T(k) = (E \oplus \mathcal{T}_{k} \otimes G_{0}^{T})^{l} \otimes \mathcal{T}_{k}
                                                 \otimes (E \oplus G_{1}^{T}).
$$

\begin{example}
Let us calculate the state transition matrix $  T(k)  $ for the system with
$  l=3 $. The above representation and idempotency in the max-plus algebra
lead us to the matrix
$$
T(k)
= \left(
    \begin{array}{cccccc}
      \tau_{1k}   & \varepsilon & \varepsilon & \varepsilon & \varepsilon
                                              & \tau_{1k} \otimes \tau_{4k} \\
      \varepsilon & \tau_{2k}   & \varepsilon & \varepsilon & \varepsilon
                            & \tau_{2k} \otimes \tau_{4k} \otimes \tau_{5k} \\
      \varepsilon & \varepsilon & \tau_{3k}   & \varepsilon & \varepsilon
          & \tau_{3k} \otimes \tau_{4k} \otimes \tau_{5k} \otimes \tau_{6k} \\
      \varepsilon & \varepsilon & \varepsilon & \varepsilon & \varepsilon
                                                                & \tau_{4k} \\
      \varepsilon & \varepsilon & \varepsilon & \varepsilon & \varepsilon
                                              & \tau_{4k} \otimes \tau_{5k} \\
      \varepsilon & \varepsilon & \varepsilon & \varepsilon & \varepsilon
                            & \tau_{4k} \otimes \tau_{5k} \otimes \tau_{6k}
    \end{array}
  \right).
$$
\end{example}

\bibliographystyle{utphys}

\bibliography{The_max-plus_algebra_approach_in_modelling_of_queueing_networks}

\end{document}